\documentclass[reqno]{amsart}
\usepackage{hyperref,amsmath,amsthm,amssymb}

\theoremstyle{plain}
\newtheorem{thm}{Theorem}
\newtheorem{lem}{Lemma}
\newtheorem{cor}{Corollary}

\renewcommand{\Re}{\mathrm{Re}}

\title{Coefficient estimates for Schwarz functions}

\author{Hitoshi Shiraishi}
\address{Hitoshi Shiraishi \newline
Department of Mathematics \newline
Kinki University \newline
Higashi-Osaka, Osaka 577-8502, Japan}
\email{step\_625@hotmail.com}

\author{Toshio Hayami}
\address{Toshio Hayami \newline
Department of Mathematics \newline
Kinki University \newline
Higashi-Osaka, Osaka 577-8502, Japan}
\email{ha\_ya\_to112@hotmail.com}

\date{}

\begin{document}
\maketitle

\begin{abstract}
Let $\mathcal{B}$ be the class of functions $w(z)$ of the form $w(z)=\sum\limits_{k=1}^{\infty}b_k z^k$ which are analytic and satisfy the condition $|w(z)|<1$ in the open unit disk $\mathbb{U}=\left\{z\in \mathbb{C}:|z|<1\right\}$.
Then we call $w(z)\in \mathcal{B}$ the Schwarz function.
In this paper, we discuss new coefficient estimates for Schwarz functions by applying the lemma due to Livingston (Proc. Amer. Math. Soc. {\bf 21}(1969), 545--552).
\end{abstract}

\

\section{Introduction}

\

Let $\mathcal{B}$ be the class of functions $w(z)$ of the form
$$
w(z)=\sum\limits_{k=1}^{\infty}b_k z^k
$$
which are analytic and satisfy the condition $|w(z)|<1$ in the open unit disk $\mathbb{U}=\left\{z\in \mathbb{C}:|z|<1\right\}$.
Also, let $\mathcal{P}$ dnote the class of functions $p(z)$ of the form
$$
p(z)=1+\sum\limits_{k=1}^{\infty}c_k z^k
$$
which are analytic and satisfy the condition $\Re(p(z))>0$ in $\mathbb{U}$.
Then we call $w(z)\in \mathcal{B}$ and $p(z)\in \mathcal{P}$ the Schwarz function and the Carath\'{e}odory function, respectively.

\

The following results are well-known for the class $\mathcal{B}$.

\

\begin{lem}[Schwarz lemma]  \label{p03lem01} \quad
If $w(z)\in \mathcal{B}$, then
$$
|w(z)|\leqq |z| \quad (z\in \mathbb{U}) \quad and \quad |b_1|\leqq 1
$$
are obtained.
In particular, $|w(z_0)|=|z_0|$ for some $z_0\in \mathbb{U}\setminus \{0\}$ or $|b_1|=1$ if and only if $w(z)=e^{i\theta}z$ for some $\theta$ $(0\leqq\theta<2\pi)$.
\end{lem}

\

By the subordination principle, we establish the following coefficient bounds.

\

\begin{lem} \label{p03lem02} \quad
If $w(z)\in \mathcal{B}$, then
$$
|b_k|\leqq 1\qquad (k=1,2,3,\ldots).
$$
Furthermore, $|b_k|=1$ for some $k$ $(k=1,2,3,\ldots)$ if and only if $w(z)=e^{i\theta}z^k$.
\end{lem}

\

Applying the Schwarz--Pick lemma (see, for example, \cite{N}), we derive the next coefficient estimate.

\

\begin{thm} \label{p03thm01} \quad
If $w(z)\in \mathcal{B}$, then
$$
|b_2|\leqq 1-|b_1|^2
$$
with equality for
$$
w(z)=\left\{
\begin{array}{ccll}
e^{i\theta}z
& & & (|b_1|=1) \\
\\
\dfrac{b_1 z+e^{i\theta}z^2}{1+e^{i\theta}\overline{b_1}z}
&=& b_1 z+(1-|b_1|^2)e^{i\theta}z^2+\ldots
& (|b_1|<1)
\end{array}
\right.
$$
for each $\theta$ $(0\leqq \theta<2\pi)$.
\end{thm}

\

In this paper, we discuss new coefficient estimates for Schwarz functions by using the following lemma due to Livingston \cite{L}.

\

\begin{lem} \label{p03lem03} \quad
If $p(z)\in \mathcal{P}$, then
$$
\left|c_s-c_t c_{s-t}\right|\leqq 2
$$
for any positive integers $s$ and $t$ $(1\leqq t<s)$.
For all $s$ and $t$, the equality is attained by the function
$$
p(z)=\dfrac{1+z}{1-z}=1+\sum\limits_{k=1}^{\infty}2z^k.
$$
\end{lem}

\

\section{Main results}

\

Our first result is contained in the following theorem by use of Lemma \ref{p03lem03}.

\

\begin{thm} \label{p03thm02} \quad
If $p(z)\in \mathcal{P}$ with $c_k=2e^{i\theta}$ $(0\leqq\theta<2\pi)$ for some $k$ $(k=1,2,3,\ldots)$, then
$$
c_{nk}=2e^{in\theta}
$$
for each $n$ $(n=1,2,3,\ldots)$.
\end{thm}

\

\begin{proof} \quad
Taking $s=2k$ and $t=k$ in Lemma \ref{p03lem03}, we see that
$$
\left|c_{2k}-c_{k}^{2}\right|=\left|c_{2k}-4e^{i2\theta}\right|\leqq 2.
$$
On the other hand, we know that $|c_{2k}|\leqq 2$, and therefore it follows that $c_{2k}=2e^{i2\theta}$.
Similarly, since
$$
\left|c_{3k}-c_k c_{2k}\right|=\left|c_{3k}-4e^{i3\theta}\right|\leqq 2\quad {\rm and}\quad |c_{3k}|\leqq 2,
$$
we have that $c_{3k}=2e^{i3\theta}$.
In the same manner, for all $n$ $(n=1,2,3,\ldots)$, we conclude that $c_{nk}=2e^{in\theta}$.
\end{proof}

\

By virtue of the above theorem, we obtain

\

\begin{cor} \label{p03cor01} \quad
If $p(z)\in \mathcal{P}$ with $c_1=2e^{i\theta}$ $(0\leqq\theta<2\pi)$, then we can declare that
$$
p(z)=\dfrac{1+e^{i\theta}z}{1-e^{i\theta}z}=1+\sum\limits_{k=1}^{\infty}2e^{ik\theta}z^k.
$$
\end{cor}

\

Moreover, applying Lemma \ref{p03lem03}, we have a new coefficient bound for Schwarz functions.

\

\begin{thm} \label{p03thm03} \quad
If $w(z)\in \mathcal{B}$, then
$$
|b_3|\leqq 1-|b_1|^3.
$$
\end{thm}

\

\begin{proof} \quad
We first note that if a function $w(z)\in\mathcal{B}$ then
$$
e^{i\theta}w(z)\in\mathcal{B}
$$
for all $\theta$ $(0\leqq\theta<2\pi)$.
Then, we know that the function $p(z)$ defined by
\begin{eqnarray} \label{p03thm03eq01}
p(z)
&=& \frac{1+e^{i\theta}w(z)}{1-e^{i\theta}w(z)} \nonumber \\
&=& 1
+ 2 e^{i\theta} b_1 z
+ 2(e^{i2\theta} b_1^2 + e^{i\theta} b_2) z^2
+ 2(e^{i3\theta} b_1^3 + 2e^{i2\theta} b_1 b_2 + e^{i\theta} b_3) z^3 \nonumber \\
& & + 2(e^{i4\theta} b_1^4 + 3e^{i3\theta} b_1^2 b_2 + 2e^{i2\theta} b_1 b_3 + e^{i2\theta} b_2^2 + e^{i\theta} b_4) z^4
+ \ldots \nonumber \\
&=& 1+c_1z+c_2z^2+c_3z^3+c_4z^4+\ldots
\qquad(z\in\mathbb{U}),
\end{eqnarray}
belongs to the class $\mathcal{P}$.
In view of Lemma \ref{p03lem03}, we obtain that
\begin{align*}
|c_3-c_1c_2|
&= |2(e^{i3\theta} b_1^3 + 2e^{i2\theta} b_1 b_2 + e^{i\theta} b_3) - 2 e^{i\theta} b_1 2(e^{i2\theta} b_1^2 + e^{i\theta} b_2)| \\
&= |2 e^{i\theta} (b_3 - e^{i2\theta} b_1^3)| \\
&\leqq 2
\end{align*}
which gives us that
$$
|b_3 - e^{i2\theta} b_1^3| \leqq 1.
$$
Thus, $b_3$ is in the region
$$
\bigcap_\theta \{ b_3 : |b_3 - e^{i2\theta} b_1^3| \leqq 1 \}
= \{ b_3 : |b_3| \leqq 1 - |b_1|^3 \}.
$$
This completes the proof of the theorem.
\end{proof}

\

The same process in the proof of Theorem \ref{p03thm03} leads us another proof of Theorem \ref{p03thm01}.

\

\begin{proof} \quad
Applying Lemma \ref{p03lem03} to the function (\ref{p03thm03eq01}) with $s=2$ and $t=1$, we deduce that
\begin{align*}
|c_2-c_1^2|
&= |2(e^{i2\theta} b_1^2 + e^{i\theta} b_2) - (2 e^{i\theta} b_1)^2| \\
&= |2 e^{i\theta} (b_2 - e^{i\theta} b_1^2)| \\
&\leqq 2.
\end{align*}
This implies that for all $\theta$ $(0\leqq\theta<2\pi)$
$$
|b_2 - e^{i\theta} b_1^2| \leqq 1
$$
which means that
$$
|b_2| \leqq 1 - |b_1|^2.
$$
\end{proof}

\

But, using the same process in the proof of Theorem \ref{p03thm03}, we have no good result for coefficients $b_4$, $b_5$ and so on.
Because, for example, applying Lemma \ref{p03lem03} to the function (\ref{p03thm03eq01}) with $s=4$ and $t=1$ to find the estimate of $b_4$, we obtain the next inequality.
\begin{eqnarray*}
|c_4-c_1c_3|
&=& |2(e^{i4\theta} b_1^4 + 3e^{i3\theta} b_1^2 b_2 + 2e^{i2\theta} b_1 b_3 + e^{i2\theta} b_2^2 + e^{i\theta} b_4) \\
& & - 4 e^{i\theta} b_1 (e^{i3\theta} b_1^3 + 2e^{i2\theta} b_1 b_2 + e^{i\theta} b_3)| \\
&=& |2 e^{i\theta} (b_4 + e^{i\theta} b_2^2 - e^{i2\theta} b_1^2 b_2 - e^{i3\theta} b_1^4)| \\
&\leqq& 2.
\end{eqnarray*}
Calculating this inequality, we have
\begin{equation} \label{p03eq01}
|b_4 + e^{i\theta} b_2^2 - e^{i2\theta} b_1^2 b_2 - e^{i3\theta} b_1^4| \leqq 1.
\end{equation}
Also, if we apply Lemma \ref{p03lem03} to the function (\ref{p03thm03eq01}) with $s=4$ and $t=2$, then we obtain another inequality as follows:
\begin{eqnarray*}
|c_4-c_2^2|
&=& |2(e^{i4\theta} b_1^4 + 3e^{i3\theta} b_1^2 b_2 + 2e^{i2\theta} b_1 b_3 + e^{i2\theta} b_2^2 + e^{i\theta} b_4) \\
& & - ( 2(e^{i2\theta} b_1^2 + e^{i\theta} b_2) )^2| \\
&=& |2 e^{i\theta} (b_4 + 2e^{i\theta} b_1 b_3 - e^{i\theta} b_2^2 - e^{i2\theta} b_1^2 b_2-e^{i3\theta} b_1^4)| \\
&\leqq& 2.
\end{eqnarray*}
Calculating this inequality, we have
\begin{equation} \label{p03eq02}
|b_4 + 2e^{i\theta} b_1 b_3 - e^{i\theta} b_2^2 - e^{i2\theta} b_1^2 b_2-e^{i3\theta} b_1^4| \leqq 1.
\end{equation}
We don't know the region to which both inequalities (\ref{p03eq01}) and (\ref{p03eq02}) point.

\

\end{document}